\documentstyle[12pt,twoside]{article}
\newtheorem{theorem}{\bf Theorem}[section]

\newtheorem{corollary}[theorem]{\bf Corollary}

\newtheorem{remark}[theorem]{\bf Remark}

\input{amssym}

\date{}
\begin{document}

\title{{\Large\bf Products of generalized derivations on rings}}

\author{{\normalsize\sc S. R.  Behresi and M. J. Mehdipour\footnote{Corresponding author}}}
\maketitle

{\footnotesize  {\bf Abstract.} In this paper, we show that if the product $(D_1D_2,
d_1d_2)$ of generalized derivations $(D_1, d_1)$ and $(D_2, d_2)$
on an algebra $A$ is a generalized derivation, then $d_1D_2$ and
$d_2D_1$
map $A$ into $\hbox{rad}(A)$. Also, for generalized derivations $(D_1, d_1)$ and $(D_2, d_2)$  on a prime ring with characteristic different from two, we give necessary and sufficient conditions under which $(D_1^2+D_1D_2,d_1^2+d_1d_2)$ is a generalized derivation as well.}
{}
\section{\normalsize\bf Introduction}

Throughout this paper $R$ denotes a ring with Jacobson radical $\hbox{rad}(R)$ and nilradical $\hbox{nil}(R)$, the intersection of all prime ideals of $R$. Let us recall that a ring $R$ is said to be
\emph{prime} if for every $x, y\in R$ the relation $xRy=\{0\}$
implies $x= 0$ or $y= 0$ and if for every $x\in R$ the relation
$xRx=\{0\}$ implies $x=0$, then $R$ is called \emph{semiprime}.
Let us also recall that an additive mapping $d: R\rightarrow R$ is
said to be a \emph{derivation} if
$$
d(xy)= d(x)y+xd(y)
$$
for all $x, y\in R$. Derivations on algebras are defined similarly.

For an algebra $A$, the mapping $F: A\rightarrow A$ is called an \emph{elementary operator} if there exist $a_i, b_i\in A$ for $i=1, ..., n$ such that
$$
F(x)=\sum_{i=1}^n a_ixb_i.
$$
These operators have an important role in the theory of operator algebras. In the case that $n=2$, $A$ has an identity $1_A$ and $b_1=a_2=1_A$ we have
$
F(x)=a_1x+xb_2
$
for all $x\in A$. Then for every $x, y\in A$
\begin{eqnarray*}
F(xy)=F(x)y+ xI_{b_2}(y),
\end{eqnarray*}
where $$I_{b_2}(t)=tb_2-b_2t$$ is the inner derivation on $A$ defined by $b_2$.
Motivated by this relation, the concept of generalized derivation on rings is introduced as follows. An additive mapping $D$ on a ring $R$ is called a \emph{generalized derivation} if there exists a derivation $d$ on $R$  such that
$$
D(xy)=D(x)y+ xd(y)
$$
for all $x, y\in R$. A generalized derivation with associated derivation $d$ is denoted by $(D, d)$. Obviously, generalized derivations include derivations and left
centralizers, an additive mapping $T$ on $R$ satisfying
$$
T(xy)= T(x)y
$$
for all $x, y\in R$. Generalized derivations  have been studied by many authors [1-4, 9-12]. For example, Hvala [11] studied the product of generalized derivations on a prime ring $R$ with $\hbox{char}(R)\neq 2$ and
characterized generalized derivations whose their product is again a
generalized derivation. Argac et al. [3] studied the product of generalized derivations $(D_1, d_1)$ and $(D_2, d_2)$ on a 2-torsion free semiprime ring $R$ and showed that $(D_1D_2, d_1d_2)$ is a generalized derivation if and only if for every $x, y\in R$,
$$
d_i(x) R D_j(y)= D_j(x) R d_i(y)=0,
$$
where $i, j\in\{1, 2\}$ and $i\neq j$; see also [13, 14].

In this paper, we continue theses investigations.
In Section 2, we investigate the product of generalized derivations on algebras and prove that if $(D_1, d_1)$, $(D_2, d_2)$ and $(D_1D_2, d_1d_2)$ are generalized
derivations on an algebra $A$, then $d_1D_2$ and $d_2D_1$ map $A$
into $\hbox{rad}(A)$. This result is an analog of Posner's first theorem for generalized derivations on an algebra. In section 3, for generalized derivations $(D_1, d_1)$ and $(D_2, d_2)$  on a prime ring with characteristic different from two, we give necessary and sufficient conditions under which $(D_1^2+D_1D_2,d_1^2+d_1d_2)$ is a generalized derivation as well.

\section{\normalsize\bf The product of generalized derivations of algebras}

The main result of this section is the following result.

\begin{theorem} Let $(D_1, d_1)$ and $(D_2, d_2)$ be generalized derivations on an algebra $A$ and $\alpha\in {\Bbb C}$.
 If $(\alpha D_1^3+ D_1D_2, \alpha d_1^3+ d_1d_2)$ is a generalized derivation on $A$, then $d_2D_1$ map $A$ into $\emph{\hbox{rad}}(A)$.
\end{theorem}
{\it Proof.} Let$(\alpha D_1^3+ D_1D_2, \alpha d_1^3+ d_1d_2)$
be generalized derivation on $A$. Then for every $x, y\in A$
\begin{eqnarray}
\alpha D_1^3(x) d_1(y)+ 3\alpha D_1(x) d_1^2(y)+ D_2(x)d_1(y)+ D_1(x)d_2(y)= 0.
\end{eqnarray}
Let $P$ be a minimal prime ideal of $A$. Since $\alpha d_1^3+ d_1d_2$ is a derivation, a similar process to that in Posner's first theorem [15] yields the fact that for every $x,
y, z\in A$
\begin{eqnarray}
d_2(x)d_1(y)d_2(z)= 0.
\end{eqnarray}
This implies that for every $x, y\in A$,
$$
d_i(x)d_j(y)\in P\quad\hbox{and}\quad d_id_j(x)\in\hbox{rad}(A)
$$
for $i, j\in\{1, 2\}$ and $i\neq j$. From this and (1) we see that
$$
d_2(D_1(x)d_2(y))\in P
$$
for all $x, y\in A$. Hence
\begin{eqnarray}
d_2(D_1(x))d_2(y)+ D_1(x)d_2^2(y)\in P
\end{eqnarray}
for all $x, y\in A$. Replacing $y$ by $d_2(y)$ in (1), and using (2), we get
\begin{eqnarray}
D_1(x)d_2^2(y)\in P
\end{eqnarray}
for all $x, y\in A$. Combining (3) and (4), gives
$$
d_2(D_1(x))d_2(y)\in P
$$
for all $x, y\in A$. It follows that $d_2(D_1(x))\in P$ for all $x\in A$. Therefore, $d_2D_1$ maps $A$ into its radical.$\hfill\square$\\

As a consequence of Theorem 2.1 we have the following result.

\begin{corollary} Let $(D_1, d_1)$ and $(D_2, d_2)$ be generalized derivations on an algebra $A$. If $(D_1D_2,d_1d_2)$ is a generalized derivation on $A$,
then $d_1D_2$ and $d_2D_1$ map $A$ into $\emph{\hbox{nil}}(A)$.
\end{corollary}
{\it Proof.} In view of proof of Theorem 2.1, we have
\begin{eqnarray}
D_2(x)d_1(y)+ D_1(x)d_2(y)= 0
\end{eqnarray}
for all $x, y\in A$. Let $P$ be a minimal prime ideal of $A$. From (5) we infer that
$$
d_1(D_2(x)d_1(y)+ D_1(x)d_2(y))\in P
$$
and
$$
D_2(x)d_1^2(y)+ D_1(x)d_2(d_1(y))\in P.
$$
Hence $d_1(D_2(x))d_1(y)\in P$ for all $x, y\in A$. Hence $d_1(D_2(x))\in P$ for all $x\in A$.$\hfill\square$

\section{\normalsize\bf The product of generalized derivations of prime rings}

We commence this section with the main result of the section.

\begin{theorem} Let $R$ be a prime
ring with $\emph{char} (R)\neq 2$ and $(D_i,d_i)$
be generalized derivations on $R$ for $i=1,2$. Then
$(D_1^2+D_1D_2,d_1^2+d_1d_2)$ is a generalized
derivation on $R$ if and only if one of the following statements holds.

\emph{(i)} $d_1\neq 0$, $d_1=-d_2$ and $D_1=-D_2$.

\emph{(ii)} $d_1=0$, $d_2\neq 0$ and $D_1=0$.

\emph{(iii)} $D_1$ and $D_2$ are right centralizer.
\end{theorem}
{\it Proof.} Let $(D_1,d_1)$ and $(D_2, d_2)$ be generalized
derivations on $R$. Then $(D_1^2+D_1D_2,d_1^2+d_1d_2)$ is a generalized derivation on $R$ if and only if
\begin{eqnarray}
2D_1(x)d_1(y)+D_2(x)d_1(y)+D_1(x)d_2(y)=0
\end{eqnarray}
for all $x,y\in R$. Hence each one of statements (i), (ii) or (iii) follows that  $(D_1^2+D_1D_2,d_1^2+d_1d_2)$ is a generalized
derivation on $R$.

Conversely, let $d_1$, $d_2$ and $d_1^2+d_1d_2$
be derivations. Then
\begin{eqnarray*}
(d_1^2+d_1d_2)(xy)&=&d_1^2(x)y+2d_1(x)d_1(y)+xd_1^2(y)\\
&+&d_1(x)d_2(y)+d_2(x)d_1(y)+xd_1d_2(y)\\
&=&(d_1^2+d_1d_2)(x)y+x(d_1^2+d_1d_2)(y)
\end{eqnarray*}
for all $x,y\in R$. Thus
\begin{eqnarray}
2d_1(x)d_1(y)+d_2(x)d_1(y)+d_1(x)d_2(y)=0
\end{eqnarray}
for all $x,y\in R$. Replace $x$ by $xd_1(z)$ in
(7). Then
\begin{eqnarray*}
2d_1(x)d_1(z)d_1(y)&+&2xd_1^2(z)d_1(y)+d_2(x)d_1(z)d_1(y)\\
&+&xd_2d_1(y)+
d_1(x)d_1(z)d_2(y)+xd_1^2(z)d_2(y)=0.
\end{eqnarray*}
From this and (7) we have
$$
2d_1(x)d_1(z)d_1(y)+d_2(x)d_1(z)d_1(y)+d_1(x)d_1(z)d_2(y)=0
$$
for all $x,y\in R$. This together with (7) shows that
$$-d_1(x)d_2(z)d_1(y)+d_1(x)d_1(z)d_2(y)=0$$
and so
$$
d_1(x)(d_1(z)d_2(y)-d_2(z)d_1(y))=0.
$$
The primeness of $R$ yields that either $d_1(x)=0$ or
$d_1(z)d_2(y)=d_2(z)d_1(y)$ for all
$x,y,z\in R$. In the both cases,
$d_1(z)d_2(y)=d_2(z)d_1(y)$ for all $y,z\in R$.
From this and (7) we infer that
$$
2d_1(z)d_1(y)+2d_2(z)d_1(y)=0.
$$
So
$$
(d_1(z)+d_2(z)) R d_1(y)=0
$$
for all $y, z\in R$. This implies that
\begin{eqnarray}
d_1=-d_2\quad\hbox{or}\quad d_1=0.
\end{eqnarray}
Now, let $D_1$ and $D_2$ both are not right centralizer. Then $d_1\neq 0$ or $d_2\neq 0$. In the case $d_1\neq 0$, we have $d_1=-d_2$ by (8). Let $d_1(y)\neq 0$ for some $y\in R$. Then from
(6) we conclude that
$$
D_1(x)d_1(y)+D_2(x)d_1(y)=0
$$
for all $x\in R$. Thus
$$
D_1(x)=-D_2(x)$$ for all $x\in R$. That is, (i) holds. Let us now consider the case $d_2\neq 0$ and $d_1=0$. Let $d_2(y)\neq 0$ for some $y\in R$.
Then by (6), $D_1(x)d_2(y)=0$ and so
$D_1(x)=0$ for all
$x\in R$.
Hence $D_1(x)=0$ for all $x\in R$. That is, (ii) holds.$\hfill\square$

\begin{remark} {\rm It is easy to see that Theorem 3.1 is not true without the assumption that $\hbox{char}(R)\neq 2$. It suffices to suppose that $D_1=D_2$ be any generalized derivation on $R$.}
\end{remark}

Now, we give some consequences of Theorem 3.1.

\begin{corollary} Let the assumptions of Theorem 3.1 be
fulfilled. Let also $(D_1^2+D_1D_2,d_1^2+d_1d_2)$ is a generalized
derivation on $R$. Then the following assertions are equivalent.

\emph{(a)} $D_1$ is a nonzero derivation;

\emph{(b)} $D_2$ is a nonzero derivation and $D_2= -D_1$;

\emph{(c)} $D_2$ and $D_1$ are nonzero derivations;

\emph{(d)} $D_1$ is a nonzero derivation and $D_1= -D_2$.
\end{corollary}
{\it Proof.} Let $D_1$ be a nonzero derivation. Then $D_1= d_1$. It follows from Theorem 3.1 (i) that $D_2$ is a nonzero derivation and $D_2= -D_1$. Therefore, (a)$\Rightarrow$(b) and (c)$\Rightarrow$(d). The implications (b)$\Rightarrow$(c) and (d)$\Rightarrow$(a) are clear.$\hfill\square$

\begin{corollary} Let $d_1$ and $d_2$ be derivations on a prime ring $R$ with $\emph{char} (R)\neq 2$. If $d_1^2+d_1d_2$ is a derivation on $R$, then $d_1=0$ or $d_1=-d_2$.
\end{corollary}

Let $\xi, \zeta\in R$. The generalized derivation $r\mapsto \zeta r+ r\xi$ on $R$ is denoted by $F_{\xi, \zeta}$ and is called an \emph{inner generalized derivation}.
In this case $(D, I_{\zeta})$ is a generalized derivation on $R$.

\begin{corollary} Let the assumptions of Theorem 3.1 be
fulfilled. Let also $(D_1^2+D_1D_2,d_1^2+d_1d_2)$ is a generalized
derivation on $R$. If $\xi, \zeta\in R$, then the following statements hold.

\emph{(i)} If $D_1= F_{\xi, \zeta}$, then $D_2= -F_{\xi, \zeta}$  or $\zeta\in\hbox{Z}(R)$.

\emph{(ii)} If $D_2= F_{\xi, \zeta}$  and $D_1\neq 0$, then $D_1= -F_{\xi, \zeta}$ or $\zeta\in\hbox{Z}(R)$.
\end{corollary}

For every $\xi\in R$, the right centralizer $r\mapsto r\xi$ on $R$ is denoted by $R_\xi$. Note that $(R_\xi, I_\xi)$ is a generalized derivation on $R$.

\begin{corollary} Let the assumptions of Theorem 3.1 be
fulfilled. Let also $(D_1^2+D_1D_2,d_1^2+d_1d_2)$ is a generalized
derivation on $R$. Then the following statements hold.

\emph{(i)} If $D_1$ is a nonzero left centralizer on $R$, then $D_2$ is a left centralizer on $R$.

\emph{(ii)} If $D_2$ is a left centralizer on $R$, then $D_1$ is a left centralizer on $R$.

\emph{(iii)} If $D_1= R_\xi$, then $D_2=R_{-\xi}$ or $\xi\in\hbox{Z}(R)$.

\emph{(iv)} If $D_2=R_\xi$ and $D_1\neq 0$, then $D_1=R_{-\xi}$ or $\xi\in\hbox{Z}(R)$.
\end{corollary}

\vspace{2mm}

 {\footnotesize
\noindent {\bf Seyed Reza  Behresi}\\
Department of Mathematics,\\ Shiraz University of Technology,\\
Shiraz
71555-313, Iran\\ e-mail: r.behresi@sutech.ac.ir\\
{\bf Mohammad Javad Mehdipour}\\
Department of Mathematics,\\ Shiraz University of Technology,\\
Shiraz
71555-313, Iran\\ e-mail: mehdipour@sutech.ac.ir\\
\end{document}